\documentclass[11pt,a4paper]{article}

\newif\ifsubmit

\usepackage[utf8]{inputenc}
\usepackage{amsmath, amsfonts, amssymb, mathtools}
\usepackage{subcaption}
\usepackage[left=2cm,right=2cm,top=2cm,bottom=2cm]{geometry}
\usepackage[bookmarks=true, pdftoolbar=true, backref=page, breaklinks]{hyperref}

\ifsubmit
	\usepackage{graphicx, tikzexternal, tikzscale}
	\graphicspath{{extfig/}}
\else
	\usepackage{tikz, tikzscale, pgfplots, pgfplotstable}
	\usetikzlibrary{plotmarks}
	\usetikzlibrary{external}
	\pgfplotsset{compat=newest} 
	\pgfplotsset{plot coordinates/math parser=false} 
	\tikzset{mark size=1}
	\tikzset{font=\footnotesize}

	\AtEndPreamble{
  	\LetLtxMacro{\oldincludegraphics}{\includegraphics}
  	\newcommand{\myincludegraphics}[2][]{\tikzsetnextfilename{#2}\oldincludegraphics[#1]{#2}}
  	\LetLtxMacro{\includegraphics}{\myincludegraphics}
	}

\fi

\usepackage{xcolor,soul, lineno}
\usepackage{algorithm, algorithmic}

\author{Shriram Srinivasan \& Anna-Karin Tornberg}
\title{Fast Ewald summation for Green's functions of Stokes flow in a half-space}
\usepackage{xspace}
\newcommand{\norm}[1]{\ensuremath{\left\lVert#1\right\rVert}}                             
\newcommand{\mb}[1]{\ensuremath{\mathbf{#1}}}                                   
\newcommand{\bs}[1]{\ensuremath{\boldsymbol{#1}}}                               
\newcommand{\mr}[1]{\ensuremath{\mathrm{#1}}}                                   
  
\newcommand{\freespace}{\emph{free-space}\xspace}
\newcommand{\halfspace}{\emph{half-space}\xspace}        
\newcommand{\etal}{\emph{et al}\xspace}

\newcommand{\stokeslet}{S}
\newcommand{\stresslet}{\mathcal{T}}
\newcommand{\rotlet}{W}
\newcommand{\harmonic}{H}
\newcommand{\biharmonic}{B}
\newcommand{\kmax}{k_{\infty}}
\newcommand{\trunc}{\mathcal{R}}

\newcommand{\xb}{\vec{x}}
\newcommand{\yb}{\vec{y}}

\newcommand{\erf}{\operatorname{erf}}
\newcommand{\erfc}{\operatorname{erfc}}

\newcommand{\Gwoexp}{A}

\begin{document}
\maketitle

\begin{abstract}
Recently, Gimbutas \etal in \cite{Gimbutas2015} derived an elegant representation for the Green's functions of Stokes flow in a half-space. We present a fast summation method for sums involving these \halfspace Green's functions (stokeslets, stresslets and rotlets) that consolidates and builds on the work  by Klinteberg \etal \cite{fsewald2017} for the corresponding \freespace Green's functions.  {The fast method is based on two main ingredients: The Ewald decomposition and subsequent use of FFTs.} The Ewald decomposition recasts the sum into a sum of two exponentially decaying series: one in real-space (short-range interactions) and one in Fourier-space (long-range interactions) with the convergence of each series  controlled by a common parameter. {The evaluation of short-range interactions is accelerated by restricting computations to neighbours within a specified distance, while the use of FFTs accelerates the computations in Fourier-space thus accelerating the overall sum.} We demonstrate  that while the  method incurs extra costs for the \halfspace in comparison to the \freespace evaluation, greater computational savings is also achieved when  compared to their respective direct sums. 
\end{abstract}
 

\section{Introduction}

The determination of the motion of particles in bounded or unbounded
flows is a central problem in microhydrodynamics. For a large class of
industrial processes like particle filtration, sedimentation or
aggregation, and deposition of pulp fibres in paper manufacturing, the
fluid inertia is negligible and the governing equations are
well-approximated by the Stokes equations \cite{KimKarilla1991}.

The system of Stokes equations is linear and can be reformulated as an
integral equation. In a boundary integral method, once the integrals
are discretized, discrete sums with fundamental solutions of Stokes
flow remain. Typically, any exterior solid boundaries or interfaces
between different fluids are discretized. Periodicity in one or more
directions is however usually built into the definition of the
fundamental solution, leading to an (infinite) summation of periodic
images in the discrete sums. Thus one can speak of 1-periodic,
2-periodic, 3-periodic or \freespace problems - indicating the
periodicity built into the evaluation of the Stokes potentials.
The finite \freespace sums for evaluating Stokes potentials are of the
form (or some variant of) \eqref{eqn:v_sum_free_space}, i.e., each summand is a convolution of
the Green’s function with a source term.

For a simple geometry like a flat plane in unbounded space, it is also
possible to avoid explicit discretization of this plane, and instead
modify the evaluation of the Stokes potentials to achieve a no-slip
condition. This requires the introduction of sources at image locations
reflected in the wall, as well as correction terms. Such an explicit
representation was first derived by Blake \cite{JBlake1971} and later in a more
elegant form by Gimbutas \etal \cite{Gimbutas2015}. 
The fast evaluation of these sums will be the focus of {this} paper, and we
start by giving some background to the problem.

If the Green's function is {the} harmonic kernel and the source term a scalar, the sum corresponds to the formula for the Coulomb potential of a system of point charges. The interest in a fast evaluation of such sums actually stemmed from this particular problem, with Ewald's investigation \cite{Ewald1921} of the 3-periodic case in 1921 now known as the Ewald summation technique. Hasimoto \cite{Hasimoto1959} then considered the 3-periodic sum of stokeslets. 
In these decompositions, some specific choices are made to turn a
conditionally convergent sum into two rapidly converging sums - one in
real space and one in Fourier space. The {computational} complexity is
however quadratic in the number of points. The survey by Deserno and
Holm \cite{Deserno1998} traces the development of fast methods based on FFTs for
acceleration of the Fourier space sum in the context of
electrostatics and molecular dynamics. An early method called Smooth
Particle Mesh Ewald (SPME) \cite{Essmann1995} was utilized by Saintillan \etal \cite{Saintillan2005} for fast evaluation of periodic stokeslet sums. 
{Later, Tornberg \etal  developed a Spectral Ewald (SE) method  for the 3-periodic sum of stokeslets \cite{Lindbo2010}, stresslets \cite{AfKlinteberg2014a} and rotlets \cite{AfKlinteberg2016rot} that is spectrally accurate and recovers the exponentially fast convergence of the Ewald sums that traditional Particle Mesh Ewald approaches cannot.}


The SE method is best suited for the 3-periodic case. Otherwise, for every direction that is not periodic, oversampling of the FFTs becomes necessary to compute the aperiodic convolution, and this increases the computational cost. The work in \cite{Lindbo2011e} illustrates the use of the SE method for a 2-periodic sum of stokeslets, while  in \cite{SE1P_Ewald_electrostatics} the SE method is adapted for 1-periodic sums in the context of electrostatics. The case of \freespace sum of stokeslets (no periodicity) is the most challenging for the SE method and it was solved recently by Klinteberg \etal \cite{fsewald2017} by combining two different ideas. The first idea is the \freespace solution of harmonic and biharmonic equations using FFT on a uniform grid by Vico \etal \cite{Vico2016} that amounts to the convolution of harmonic/biharmonic (radial) kernels with source terms by FFT on a uniform grid, and the second idea is that the stokeslet, stresslet and rotlet kernels, though not radial, can be expressed as a linear combination of differential operations on the harmonic or biharmonic kernels. 

A popular method ideally suited for \freespace problems is the Fast Multipole Method (FMM) \cite{Greengard1987}. In contrast to the SE method, it is best suited for problems with no periodicity. The FMM has been used successfully for harmonic and biharmonic  kernels and stokeslets \cite{Rodin2000, Duraiswami2006, Tornberg2008}. However, the SE method still compares favourably with FMM for \freespace \cite{fsewald2017}. While this comparison was based on a uniform distribution of source points, for increasingly non-uniform distribution of sources, the adaptivity of the FMM will eventually come into play and be a decisive advantage. Irrespective of this  observation, one of the merits of the SE method is its versatility, having been shown to work for 3-periodic, 2-periodic and 
\freespace Stokes flow problems. 

In many studies of sedimentation, it is natural to consider a
\halfspace ($\mathbb{R}^2 \times \mathbb{R^+}$) domain bounded by a plane wall at the bottom.
As mentioned above, this can be achieved without explicitly
discretizing the wall, instead modifying the discrete sums to be
evaluated \cite{Gimbutas2015}.

This work deals with the following specific problem: that of fast
summation of a large number of convolutions of the Green’s function
for a \halfspace with source terms. We apply Ewald decomposition to
the modified formulae, and adapt the SE method to this case.  The
structure of the Green’s function for the \halfspace is however more
complicated than other kernels investigated previously.


The outline of this paper is as follows: We introduce our notation and the \freespace problem in Section~\ref{sec:freespace}. The \halfspace problem and its Green's function are stated in Section~\ref{sec:halfspace}. A brief explanation of the Ewald summation technique as it applies here is given in Section~\ref{sec:Ewald} before we describe the Spectral Ewald Method for the case of the \halfspace in Section~\ref{sec:hse}. Some error estimates are discussed in Section~\ref{sec:truncation} before we present computational results in Section~\ref{sec:results}.

\section{The free-space problem}
\label{sec:freespace}
The three fundamental solutions for Stokes flow are the so-called stokeslet, stresslet and rotlet singularities; they are tensors that have an explicit representation as follows:
\begin{subequations}
\begin{align}
8\pi \mb{S} = &  \dfrac{\bs{I}}{r} + \dfrac{\bs{r} \otimes \bs{r} }{r^3}, \\
8\pi \mathcal{T} = & -6 \dfrac{\bs{r} \otimes \bs{r} \otimes \bs{r}}{r^5}, \\
4\pi \mb{W} =& \dfrac{\mathcal{E}\bs{r}}{r^3},
\end{align}
\end{subequations}
with the short hand notation 
$\norm{\mb{r}} =: r$ and $\mb{x}, \; \mb{y} \in \mathbb{R}^3$. Here $\otimes$ is the standard tensor product for vectors in $\mathbb{R}^3$, $\bs{I}$ is the second order identity tensor in $\mathbb{R}^3$ and $\mathcal{E}$ is the alternating third-order tensor whose representation in the natural basis of $\mathbb{R}^3$ coincides with that of Levi-Cevita's symbol. 

In the Einstein summation convention or index notation, the direct representations above reduce to
\begin{subequations}
\begin{align}
8\pi S_{ij} = &  \dfrac{\delta_{ij}}{r} + \dfrac{r_i r_j}{r^3}, \\
8\pi T_{ijk} = & -6 \dfrac{r_i r_j r_k}{r^5}, \\
4\pi W_{ij} = & \dfrac{\epsilon_{ijk}r_k}{r^3},
\end{align}
\end{subequations}
where $\delta_{ij}, \epsilon_{ijk}$ are the Kronecker delta and Levi-Cevita's symbol respectively. Wherever possible we shall use the direct notation for economy and provide the index notation as explanation.

Given that there are $N$ point forces of intensity $8\pi\mb{f}^{(m)}$ at locations $\mb{y}^{(m)}, \; m = 1\dotsc N$,  the  stokeslet induces a velocity field at $\mb{x} \in \mathbb{R}^3$. Similarly, if there are $N$  ``double forces" of intensity $8\pi\mb{g}^{(m)}$ and orientation $\mb{q}^{(m)}$, the stresslet and rotlet tensors also induce a velocity as follows: 

\begin{subequations}
\begin{align}
\mb{u}(\mb{x}) = & \sum_{\substack{m=1}}^{N}\mb{S}^{(m)} \mb{f}^{(m)}, \\
\mb{u}(\mb{x})  = & \sum_{\substack{m=1}}^{N} \mathcal{T}^{(m)}(\mb{g}^{(m)} \otimes \mb{q}^{(m)}), \\
\mb{u}(\mb{x})  = & \sum_{\substack{m=1}}^{N} 2\mb{W}^{(m)}(\mb{g}^{(m)} \times \mb{q}^{(m)}), 
\end{align}
\end{subequations}
where $\mb{S}^{(m)}= \mb{S}(\bs{r}^{(m)})$, and $\mathcal{T}^{(m)}$ and $\mb{W}^{(m)}$ are defined similar to $\mb{S}^{(m)}$ with the short hand notation 
$\mb{x}- \mb{y}^{(m)} =: \bs{r}^{(m)}$.  The formula for the rotlet  can also be written in terms of the vector cross product by replacing the skew-symmetric tensor with the axial vector as has been done in \cite{fsewald2017}, but we do not do so in order to stay consistent with the notation in Gimbutas \etal \cite{Gimbutas2015}. The case where the evaluation point $\mb{x} = \mb{y}^{(p)}, p = 1, 2, \dotsc N, p\neq m$ will be the one considered here on as it is of most interest and occurs in boundary integral methods and  potential methods for Stokes equations.

In index notation, the velocity field induced in each case is written as
\begin{subequations}
\begin{align}
u_i = & \sum_{\substack{m=1}}^{N}S^{(m)}_{ij} {f}^{(m)}_j, \\
u_i = & \sum_{\substack{m=1}}^{N}T^{(m)}_{ijk} {g}^{(m)}_j {q}^{(m)}_k, \\
u_i = & \sum_{\substack{m=1}}^{N}2W^{(m)}_{ij} \epsilon_{jkl} {g}^{(m)}_k {q}^{(m)}_l. 
\end{align}
\label{eqn:v_sum_free_space}
\end{subequations}

The result of summation over all $N$ terms above in \eqref{eqn:v_sum_free_space} gives the \freespace velocity corresponding to a collection of stokeslets, stresslets or rotlets. Note that the operation in each case is a convolution\footnote{This is clear when the source $\mb{f}^{(m)} = \mb{f}(\mb{y}^{(m)})$ is written in full as  $\mb{u}(\mb{x}) = \sum_{\substack{m=1}}^{N}\mb{S}(\mb{x}- \mb{y}^{(m)}) \mb{f}^{(m)}\delta(\mb{y}^{(m)})$ } between the kernel and the source term. In \cite{fsewald2017}, a fast SE method is proposed for this evaluation. 

\section{The half-space problem}
\label{sec:halfspace}

A natural extension of this question would be to ask  if there exists an  explicit representation for the velocity field due to discrete singularities in the \halfspace, and secondly, if  there are fast methods to evaluate it.
This paper deals with the latter question. Indeed, there are closed form expressions for the \halfspace, those given by Blake \cite{JBlake1971} and Gimbutas \etal \cite{Gimbutas2015} which use the method of images to derive  an appropriate expression. While the older representation of Blake requires evaluation of multiple harmonic and dipole fields, the recent work by Gimbutas \etal provides a very elegant representation using correction terms expressed in terms of a single harmonic potential. Their key idea was to invoke the Papkovitch-Neuber \cite{papkovich1932solution, neuber1934neuer} representation formula for constructing a divergence-free velocity field that satisfies the Stokes equation using a harmonic potential. The \freespace formula itself does not satisfy the  boundary conditions at the wall, but the combination with its image ensures the tangential component satisfies the boundary condition. The Papkovitch-Neuber correction term is added to adjust the normal component of the velocity at the wall. The result is the formula \eqref{eqn:v_sum_halfspace} which is divergence-free, satisfies Stokes equations, and the no-slip boundary conditions at the wall. In this paper, we show how the SE method illustrated in \cite{fsewald2017} can be applied to this representation formula for the \halfspace \eqref{eqn:v_sum_halfspace} to yield a fast method of evaluation.

\begin{figure}[ht]
\centering
\includegraphics[scale=0.5]{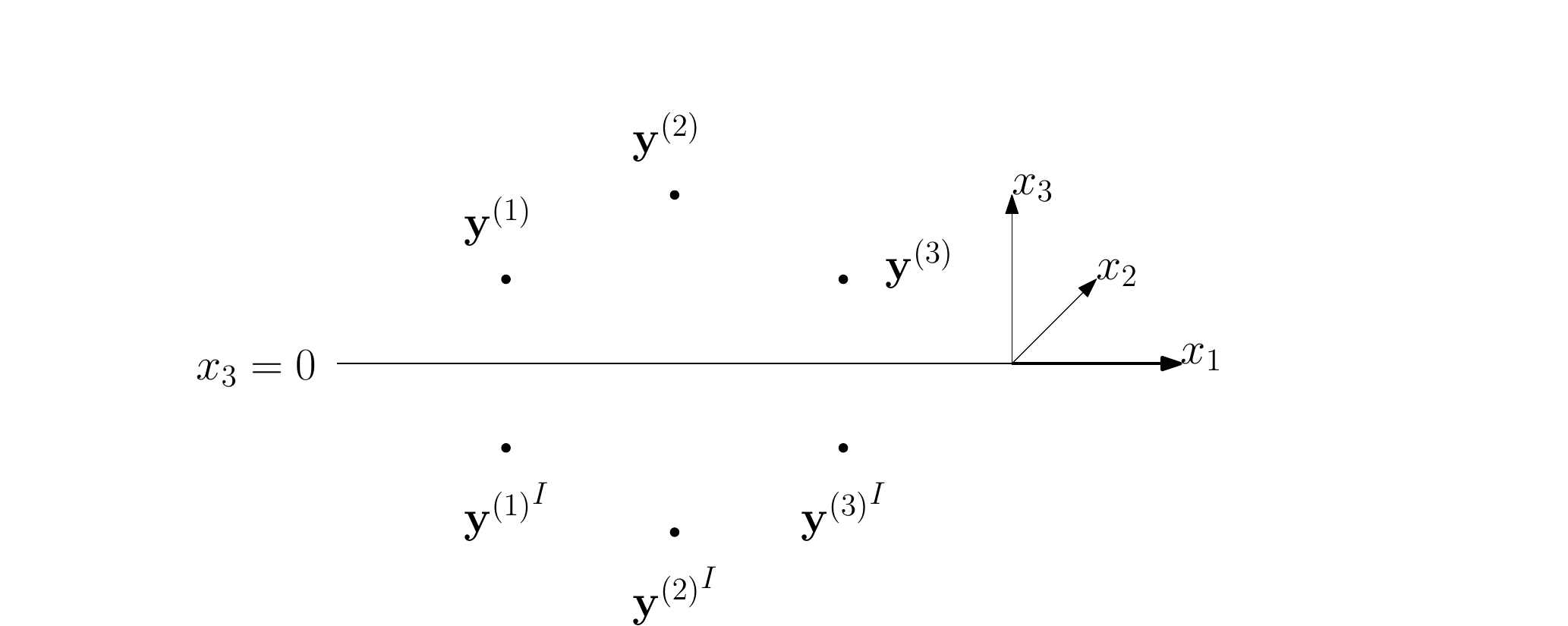}
\caption{The source locations and image locations along with the infinite wall forming the half-space}
\label{fig:halfspace}
\end{figure}

First we set out our notation in order to restate the formulae for the \halfspace. A system of Cartesian coordinates is set up such that the $x_3 = 0$ plane coincides  with the wall.

Let $\mb{y}^I := (y_1, y_2, -y_3)^T$ denote the reflection of $\mb{y} := (y_1, y_2, y_3)^T$ about the  $x_3$ plane (see Figure~\ref{fig:halfspace}). Furthermore, let  $\mb{\tilde{r}}^{(m)} := \mb{x} - {\mb{y}^{(m)}}^{I}, \; \mb{\tilde{S}}^{(m)} := \mb{S}(\mb{\tilde{r}}^{(m)})$

Then the representation for the \halfspace in each case involves a harmonic function associated with each mirror image as follows:
\begin{subequations}
\begin{align}
\mb{u}(\mb{x}) = & \sum_{\substack{m=1}}^{N}\mb{S}^{(m)} \mb{f}^{(m)} - \mb{\tilde{S}}^{(m)} {\mb{f}^{(m)}}^I -x_3\nabla \phi^{(m)}_{S} + (0, 0, \phi^{(m)}_{S}  )^T, \\
\mb{u}(\mb{x}) = & \sum_{\substack{m=1}}^{N}\mathcal{T}^{(m)}(\mb{g}^{(m)} \otimes \mb{q}^{(m)}) - 
\mathcal{\tilde{T}}^{(m)}({\mb{g}^{(m)}}^I \otimes {\mb{q}^{(m)}}^I) - x_3\nabla \phi^{(m)}_{\stresslet} + (0, 0, \phi^{(m)}_{\stresslet}  )^T, \\
\mb{u}(\mb{x}) = & \sum_{\substack{m=1}}^{N}2\mb{W}^{(m)} (\mb{g}^{(m)} \times \mb{q}^{(m)}) -
2\mb{\tilde{W}}^{(m)} ( {\mb{g}^{(m)}}^I \times {\mb{q}^{(m)}}^I) - x_3\nabla \phi^{(m)}_{W} + (0, 0, \phi^{(m)}_{W}  )^T.
\end{align}
\label{eqn:v_sum_halfspace}
\end{subequations}
The harmonic functions $\phi_S, \phi_{\stresslet}$ and $\phi_W$  presented in \cite{Gimbutas2015} are written in terms of potentials due to point sources, dipoles and quadrupoles. Here we instead express them as gradients of the harmonic potential. To explain, suppose the harmonic potential $G$ of a unit charge  $4\pi G:= \dfrac{1}{r}$, we can write the potential of a dipole with orientation $\mb{a}$ as $G^D := -\nabla G \cdot \mb{a}$ and that of a quadrupole with intensity $\mb{b}$ and orientation $\mb{a}$ as $G^Q := \nabla \nabla G \mb{a} \cdot \mb{b}$. Thus, using these relations the expressions for $\phi_S, \phi_T$ and $\phi_W$ recorded in \cite{Gimbutas2015} are restated as

 \begin{subequations}
 \begin{align}
 \phi_S^{(m)} =& 2\left( -f^{(m)}_3 G^{(m)} - y^{(m)}_3 \nabla G^{(m)} \cdot {\mb{f}^{(m)}}^I \right), \\
 \phi_T^{(m)} =& 2\left( 2 \left( {\mb{g}^{(m)}}^I \cdot {\mb{q}^{(m)}}^I \right) \left(\nabla G^{(m)} \cdot \mb{e}_3 \right) +  y^{(m)}_3 \left( \nabla \nabla G^{(m)} {\mb{g}^{(m)}}^I \cdot {\mb{q}^{(m)}}^I \right)  \right), \\
 \phi_W^{(m)} =&  2\left(  \nabla G^{(m)} \cdot 2 \left( q^{(m)}_3{\mb{g}^{(m)}}^I  -  g^{(m)}_3{\mb{q}^{(m)}}^I \right) \right),
 \end{align}
 \label{eqn:phi}
 \end{subequations}
 where $G^{(m)} = G(\tilde{r}^{(m)})$.

 Note that the formulae \eqref{eqn:v_sum_halfspace} are not translation-invariant, and it is essential that the origin of the coordinate system be located on the wall or boundary of the \halfspace.
 
\section{The Ewald decomposition} 
\label{sec:Ewald}

Here we quickly summarize the motivation behind the Ewald decomposition.
The idea is to introduce a scalar (Ewald) parameter $\xi$ and split the fundamental solution or Green's function (for the stokeslet, say) into

$$\mb{S}(\mb{x} - \mb{y}^{(m)}) = \mb{S}_{R}(\mb{x} - \mb{y}^{(m)}, \xi) +  \mb{S}_{F}(\mb{x} - \mb{y}^{(m)}, \xi) = \mb{S}_{R}(\mb{x} - \mb{y}^{(m)}, \xi)  +  
\mathcal{F}^{-1}(\widehat{\mb{S}_{F}}(\mb{k}, \xi)),$$
where $\widehat{\mb{S}_{F}}$ is the Fourier transform of $\mb{S}_{F}$ and $\mathcal{F}^{-1}$ indicates the inverse Fourier transform operator (IFT).
The formula for $\mb{u}(\mb{x})$  then becomes
$$\mb{u}(\mb{x}) = \sum_{\substack{m=1}}^{N}\mb{S}_{R}(\mb{x} - \mb{y}^{(m)}, \xi)\mb{f}^{(m)}  + \frac{1}{(2\pi)^3} \int_{\mathbb R^3} 
    \widehat{\mb{S}_{F}}(\mb{k}, \xi)) \sum_{m=1}^N \mb{f}^{(m)} e^{i \mb{k}\cdot (\mb{x} - \mb{y}^{(m)}) } \mathrm{d}\mb{k}.$$
The first term represents local interactions and it will be seen below that it decays exponentially fast, and can hence be truncated. Thus it is evaluated directly by summation over all sources located within a chosen distance $r_c$ of the target $\mb{x}$. The second term is the IFT integral and if the integrand were smooth  and compactly supported, the  integral could be approximated to spectral accuracy with a trapezoidal rule in each coordinate direction, allowing for the use of FFTs for the evaluation. However, we will find that all the kernels and correction functions relevant to this present work will have a factor (like $\widehat{B}(k)$ in the expression for $\widehat{\mb{S}_{F}}(\mb{k}, \xi)$ recorded below)  which makes them singular at $k=0$. The method to circumvent this singularity in our quadrature is one of the key points and will be discussed in detail subsequently.

 If the target location $\mb{x}$ coincides with a source location, then the so-called \emph{self-interaction} term must also be accounted for and removed. This is evaluated as the limit of $\mb{S}_{R} - \mb{S}$ when $r^{(m)} = \norm{\mb{x} - \mb{y}^{(m)}} \rightarrow 0$. 

For purpose of illustration, we explicitly show the Ewald decomposition of the stokeslet derived in  \cite{Hasimoto1959}.

$$ \mb{S}_{R}(\mb{x}, \xi) =  2\left( \frac{\xi  e^{-(\xi r^{(m)})^2 }  }{\sqrt{\pi}} +
      \frac{ \erfc{ ( \xi r^{(m)} ) } } { 2 r^{(m)} } \right) \left( \bs{I} + 
      \frac{ \bs{r}^{(m)} }{r^{(m)}} \otimes \frac{ \bs{r}^{(m)} }{r^{(m)}} \right) -
    \frac{ 4 \xi e^{-(\xi r^{(m)})^2} }{\sqrt{\pi}}    \bs{I}.  $$

$$\widehat{\mb{S}_F}(\mb{k}, \xi) =  \mb{A}(\mb{k}, \xi)e^{-k^2/4\xi^2} \hat{B}(k), $$ 
where $\mb{A}(\mb{k}, \xi) := -\left(k^2\bs{I} - \mb{k} \otimes \mb{k}
  \right) \left( 1 + \frac{k^2}{4\xi^2} \right), \;\; \mb{k} = (k_1, k_2, k_3)^T$, $k := \norm{\mb{k}}, \hat{B}(k) := -\frac{8\pi}{k^4}$.

$$ \mb{S}_{\mathrm{self}}(\xi) = \left[ -\frac{4\xi}{\sqrt{\pi}}\bs{I}\right].$$   

The corresponding expressions for the stresslet and rotlet are tabulated for reference in the appendix in equations \eqref{eqn:kernel_real_parts} and \eqref{eqn:kernel_fourier_space_part}.

The computational complexity of the direct sum in  \eqref{eqn:v_sum_free_space} is $O(N^2)$ for $\mb{x} = \mb{y}^{(p)}, p = 1, 2, \dotsc N, p\neq m$, and the Ewald sum by itself does not reduce the complexity. The SE method keeps the real-space sum at $O(N)$ by a specific choice/scaling of parameters when scaling up the 
system, and this combined with the use of FFTs reduces the over-all complexity (see \cite{fsewald2017} for details) to $O(N\log N)$.

One would like to use the SE method for \eqref{eqn:v_sum_halfspace} also to accelerate the process. While at first glance, the expressions involved in the sum may appear to be too cumbersome and complicated, a careful perusal will reveal that there is a structure to the formula; the summand consists of a linear combination of a stokeslet/stresslet/rotlet and its image about the $ x_3 = 0$ plane along with two other terms involving the respective  harmonic functions, which could be viewed as corrections to satisfy the boundary conditions.
The Ewald decomposition of the first two terms on the right hand side of \eqref{eqn:v_sum_halfspace} follows directly from \cite{fsewald2017}, so it only remains to derive an appropriate Ewald decomposition for the two correction terms involving the harmonic function and its gradient. We recognize that the Ewald decomposition of $G, \nabla G, \nabla \nabla G$ and $\nabla\nabla\nabla G$ is required to complete this task.

Turning to it, a lengthy (but straight-forward) calculation finds the real space components (denoted by the subscript $R$) to be
\begin{subequations}
\begin{align}
G_R =& \dfrac{1}{r}-f(r), \\
(\nabla G_R)_i =& -\left[ \dfrac{1}{r^3} + \dfrac{f'(r)}{r} \right]r_i, \\
(\nabla\nabla G_R)_{ij} =& -\left[ \dfrac{1}{r^3} + \dfrac{f'(r)}{r} \right]\delta_{ij} + \left[ \dfrac{3}{r^5} + \dfrac{f'(r)}{r^3} - \dfrac{f''(r)}{r^2} \right]r_ir_j, \\
(\nabla\nabla\nabla G_R)_{ijk} =&  \left[ \dfrac{3}{r^5} + \dfrac{f'(r)}{r^3} - \dfrac{f''(r)}{r^2} \right] (\delta_{ij}r_k + \delta_{ik}r_j + \delta_{jk}r_i) \\ \notag 
  & + \left[ -\dfrac{15}{r^7} -\dfrac{3f'(r)}{r^5} +\dfrac{3f''(r)}{r^4} -\dfrac{f'''(r)}{r^3} \right] r_ir_jr_k,
\end{align}
\label{eqn:G_R}
\end{subequations}
where, $r_i$ is the $i^{\textrm{th}}$ component of the vector joining the source-location to the target-location, and  $f(r) := \dfrac{\erf(\xi r)}{r}$ for convenience.

The components evaluated in Fourier-space through their Fourier Transforms and denoted by the subscript $F$ are
\begin{subequations}
\begin{align}
& \widehat{G_F}(\mb{k}, \xi) = \widehat{H}(k) e^{-\frac{k^2}{4\xi^2}}, \\
& (\widehat{\nabla G_F})_m = i k_m \widehat{G_F}(\mb{k}, \xi) = i k_m \widehat{H}(k)e^{-\frac{k^2}{4\xi^2}} =  i \hat{k}_m k\widehat{H}(k)e^{-\frac{k^2}{4\xi^2}},  \\
& (\widehat{\nabla \nabla G_F})_{mn} = -k_m k_n \widehat{G_F}(\mb{k}, \xi) = - k_m k_n \widehat{H}(k) e^{-\frac{k^2}{4\xi^2}} = - \hat{k}_m \hat{k}_n k^2\widehat{H}(k)e^{-\frac{k^2}{4\xi^2}},\\
& (\widehat{\nabla \nabla \nabla G_F})_{mnp} = -i k_m k_n k_p \widehat{G_F}(\mb{k}, \xi) = -i k_m k_n k_p \widehat{H}(k)e^{-\frac{k^2}{4\xi^2}} = -i\hat{k}_m \hat{k}_n \hat{k}_p k^3\widehat{H}(k)e^{-\frac{k^2}{4\xi^2}},
\end{align}
\label{eqn:G_F}
\end{subequations}
where $\mb{k} = (k_1, k_2, k_3)^T$, $k := \norm{\mb{k}}$, $\hat{k}_i := \dfrac{k_i}{k}$ and $\widehat{H}(k) := \frac{4\pi}{k^2}$. {The reason for this notation will become clear later.}

These correction terms are centered at the image locations but the evaluation point $\mb{x}$ is never an image location hence there is no need to account for the so-called self-interaction term. 

\section{The half-space Spectral Ewald method}
\label{sec:hse}

The \halfspace formula described in Section~\ref{sec:halfspace} essentially transforms the problem into a new \freespace problem, albeit with extra terms.  Hence the computational framework used is the same as described in \cite{fsewald2017}, following the recent idea introduced by Vico \etal. \cite{Vico2016} to solve free space problems by FFTs on uniform grids. Thus we keep our explanation of the method very general and brief, directing readers to \cite{fsewald2017} for complete details while we highlight only the differences from it and finer points of note.
A simple but useful observation is that for computations, the first two terms (kernel and its image) can be combined into a single term by absorbing the negative sign into the image source vectors, so that one can think of a kernel with $2N$ sources.

\subsection{The real space component} 

The real space components for the stresslet and rotlet derived in \cite{fsewald2017, AfKlinteberg2016rot} are reproduced in  Appendix~\ref{sec:appendix}, while that of the respective correction functions $\phi_S, \phi_{\stresslet}, \phi_W$ and their gradients can be written down directly by utilizing \eqref{eqn:G_R} in \eqref{eqn:phi}. 

A cell-list of nearest neighbours within the cut-off radius $r_c$ is prepared for each target, the difference with \cite{fsewald2017} being that all target locations were also source locations in a \freespace evaluation, but for the case of a \halfspace, target locations are  on one side of the wall and also serve as sources, but  their reflections act as sources only and are located on the other side of the wall.  Other than that, the procedure to evaluate the local interactions is as before. The choice of the cut-off radius $r_c$ is made from the desired error-level using the truncation error estimates discussed in Section ~\ref{sec:truncation}.

\subsection{The Fourier-space component}
The Fourier-space component is evaluated through the integral which always has the form
{\begin{equation}
\mb{u}_F(\mb{x}, \xi) = \frac{1}{(2\pi)^3} \int_{\mathbb R^3} e^{i \mb{k}\cdot \mb{x}}
    \mb{A}(\mb{k}, \xi) e^{-\frac{k^2}{4\xi^2}}\hat{\mathcal{K}}(k)  \sum_{m=1}^N \mb{c}^{(m)} e^{-i \mb{k}\cdot \mb{y}^{(m)} } \mathrm{d}\mb{k},
    \label{eqn:uF}
\end{equation}}
where the quantities $\mb{A}(\mb{k}, \xi),\hat{\mathcal{K}}(k)$ depend on the choice of kernel (see Appendix~\ref{sec:appendix}, $\hat{\mathcal{K}}$ is either $\hat{B}(k)$ or $\hat{H}(k)$),  while the  quantity $\mb{c}^{(m)}$ is a scalar, vector or tensor that depends on the $m^{\mr{th}}$ source term and location.
For ease of explanation we write out in full the expression that emerges for the stokeslet from \eqref{eqn:v_sum_halfspace},
\begin{align}
\mb{u}_F(\mb{x}, \xi) = & \frac{1}{(2\pi)^3} \int_{\mathbb R^3} e^{i \mb{k}\cdot \mb{x}}
    {A^S}(\mb{k}, \xi) e^{-\frac{k^2}{4\xi^2}}\hat{B}(k) \sum_{m=1}^N \left\{ \mb{f}^{(m)} e^{-i \mb{k}\cdot \mb{y}^{(m)} } -  {\mb{f}^{(m)}}^I e^{-i \mb{k}\cdot {\mb{y}^{(m)}}^I } \right\} \mathrm{d}\mb{k} \notag \\
    & -\frac{x_3}{(2\pi)^3} \int_{\mathbb R^3} e^{i \mb{k}\cdot \mb{x}} \sum_{m=1}^N \widehat{\nabla\phi^{(m)}_S}  e^{-i \mb{k}\cdot {\mb{y}^{(m)}}^I }\mathrm{d}\mb{k}  \notag \\
    &   + \frac{1}{(2\pi)^3} \int_{\mathbb R^3} e^{i \mb{k}\cdot \mb{x}} \left(0, 0, 1 \right)^T \sum_{m=1}^N \widehat{\phi^{(m)}_S} e^{-i \mb{k}\cdot {\mb{y}^{(m)}}^I} 
    \mathrm{d}\mb{k}.
    \label{eqn:uF_stokeslet}
\end{align}

Substituting from \eqref{eqn:phi} and \eqref{eqn:G_F}, it expands to

\begin{align}
\mb{u}_F(\mb{x}, \xi) = & \frac{1}{(2\pi)^3} \int_{\mathbb R^3} e^{i \mb{k}\cdot \mb{x}}
    {A^S}(\mb{k}, \xi) e^{-\frac{k^2}{4\xi^2}}\hat{B}(k)   \sum_{m=1}^{2N} \left\{ \mb{F}^{(m)} \right\} e^{-i \mb{k}\cdot \mb{z}^{(m)} } \mathrm{d}\mb{k} \notag \\
     & -\frac{x_3}{(2\pi)^3} \int_{\mathbb R^3} e^{i \mb{k}\cdot \mb{x}} \left(\mb{k}\right)  e^{-\frac{k^2}{4\xi^2}}  i\hat{H}(k)
    \sum_{m=1}^N \left\{f_3^{(m)} \right\} e^{-i \mb{k}\cdot {\mb{y}^{(m)}}^I }\mathrm{d}\mb{k} \notag \\
   & -\frac{x_3}{(2\pi)^3} \int_{\mathbb R^3}e^{i \mb{k}\cdot \mb{x}} (\mb{k}\otimes\mb{k})  e^{-\frac{k^2}{4\xi^2}}  i^2\hat{H}(k) \sum_{m=1}^N \left\{y_3^{(m)}{\mb{f}^{(m)}}^I\right\} e^{-i \mb{k}\cdot {\mb{y}^{(m)}}^I }\mathrm{d}\mb{k} \notag \\
    &   + \frac{\left(0, 0, 1\right)^T}{(2\pi)^3} \int_{\mathbb R^3} e^{i \mb{k}\cdot \mb{x}}  (1) e^{-\frac{k^2}{4\xi^2}} \hat{H}(k) \sum_{m=1}^N \left\{f_3^{(m)}\right\} e^{-i \mb{k}\cdot {\mb{y}^{(m)}}^I }\mathrm{d}\mb{k} \notag \\
    & + \frac{\left(0, 0, 1\right)^T}{(2\pi)^3} \int_{\mathbb R^3} e^{i \mb{k}\cdot \mb{x}}  
    \left( \mb{k} \right) e^{-\frac{k^2}{4\xi^2}}   i\hat{H}(k)\cdot \sum_{m=1}^N 
    \left\{ y_3^{(m)}{\mb{f}^{(m)}}^I \right\} e^{-i \mb{k}\cdot {\mb{y}^{(m)}}^I }\mathrm{d}\mb{k},
    \label{eqn:uF_stokeslet_expanded}
\end{align}
where we have rewritten the first integral as a sum over $2N$ terms by setting 
{$$\mb{F}^{(m)} = \begin{cases}
\mb{f}^{(m)} \; \forall \; 1\leq m \leq N. \\
-{\mb{f}^{(m-N)}}^I \; \forall \; N+1 \leq m \leq 2N.
\end{cases}, \quad \mb{z}^{(m)} = \begin{cases}
\mb{y}^{(m)} \; \forall \; 1\leq m \leq N. \\
{\mb{y}^{(m-N)}}^I \; \forall \; N+1 \leq m \leq 2N.
\end{cases}$$}

While the equation above has been  written out in full for the stokeslet, the method that will be discussed carries over without modification for the stresslet and rotlet as well since the  expanded formula in that case too has the same structure and number of terms. 
The explicit presentation of the five integrals that need evaluation on the right hand side of \eqref{eqn:uF_stokeslet_expanded}   serves a dual purpose: \\
 (1) It shows that the integrals arising from the correction terms have the same form as that of the integral arising from the stokeslet in \freespace. For a single integral of the type considered here, Klinteberg \etal \cite{fsewald2017} in Section~5.1 illustrate and justify the sequence of operations followed. In the case of the \halfspace however, we have five integrals of that type, and we shall now explain how to combine them together in the evaluation while following the same procedure. \\
 (2) It also makes clear that all the integrands   have the factor $\widehat{B}(k)$ or $\widehat{H}(k)$ which make it singular at $k=0$. 

 {To circumvent this issue, we will continue to follow \cite{fsewald2017} and \cite{Vico2016}  and introduce modified Green’s functions where the Fourier transforms of these functions have no singularity at $k = 0$. The truncated Green's functions are denoted by the superscript $\mathcal{R}$ that stands for the radius of the support in real space, and their Fourier transforms are given by}
\begin{align}
& \widehat H^\mathcal{R}(k) = 8 \pi 
  \left( \frac{\sin(\mathcal{R} k/2)}{k}\right)^2, \\
 & \widehat B^\mathcal{R}(k) = 4 \pi \frac{
    (2-\mathcal{R}^2k^2)\cos(\mathcal{R} k) + 2 \mathcal{R}^2 k \sin(\mathcal{R}^2 k) - 2
  }{k^4}.
\end{align}

{Using these truncated Green's functions will yield exactly the same
result for the harmonic/biharmonic equation as the original ones, as
long as the right hand side has compact support within the solution
domain, and $\mathcal{R}$ is chosen sufficiently large.
Specifically, if the solution domain is such that the largest point to
point distance is $\mathcal{R}_{\mr{max}}$, then we need $\mathcal{R} \ge \mathcal{R}_{\mr{max}}$.} 

We assume that the physical domain is a cube of size $L$ that encloses all the sources, the computational domain is a cuboid of dimensions $L \times L \times 2L$ that  encloses both the sources and their images. 

{In the Ewald decomposition, sources are convolved with Gaussians or modified Gaussians to form the right hand side that defines the Fourier space problem. 
Hence, the assumption above of a compactly supported right hand side is violated.  
In the actual discretization, we interpolate point sources to the grid using a window function. 
In the SE method, this window function is a (suitably scaled, hence not the same) Gaussian, as will be explicitly defined in Algorithm~\ref{alg:SE_halfspace}.
The domain length $L$ will be extended to $\tilde{L}$ to accommodate the support around the source locations.  
With the parameter choices that we will soon detail, this extension is sufficient and further extension will not reduce the total error. 
A more detailed discussion of this issue can be found in section 5.3 of \cite{fsewald2017} .}
%
%
The resulting computational domain is discretized by an  equi-spaced grid with spacing $h$ containing $\tilde{M} \times \tilde{M} \times 2\tilde{M} = 2{\tilde{M}}^3$ points.
Note that such a grid induces, in its $k$-space counterpart, a spacing $\Delta k = 2\pi/ \tilde{L} $.

 Fundamentally, we are performing an aperiodic convolution in Fourier space, and hence need a twice over-sampled representation of  $\widehat H^\mathcal{R}(k)$, i.e., a representation with $k$-space resolution of $\Delta k/2$ . It has been shown earlier by Vico \etal \cite{Vico2016} that the terms $\widehat H^\mathcal{R}(k), \widehat B^\mathcal{R}(k)$ can be evaluated knowing only the value of $\mathcal{R}$, which itself is determined by the size or extent of the domain. Thus for computational efficiency, we precompute $\widehat H^\mathcal{R}(k)$, and $\widehat B^\mathcal{R}(k)$ for the stokeslet or stresslet on a grid with $16{\tilde{M}}^3$ points. We cannot compute this directly by starting from  values on the physical grid since it is only the Fourier transform that is known analytically. This computation is thus carried out as follows:
\begin{enumerate}
\item {Evaluate $\widehat H^\mathcal{R}(k)$ and  $\widehat B^\mathcal{R}(k)$ on a grid of spacing $\Delta k/s_f$ (or $2(s_f \tilde{M})^3$ points) where the truncation radius for the
domain $\mathcal{R} = \sqrt{6}\tilde{L}$ and the oversampling factor $s_f \ge 1 + \sqrt{6}$ is chosen as small as possible such that $s_f \tilde{M}$ is an even integer.}
\item Compute the 3D-IFFT and truncate to get $H^\mathcal{R}$ on a grid of $16{\tilde{M}}^3$ points.
\item Compute the 3D FFT now to get back $\widehat H^\mathcal{R}(k)$ on a grid with spacing $\Delta k/2$, that is a twice oversampled representation.
\end{enumerate}
This set of values is now used in the algorithm. Note that this is different from simply sampling $\widehat H^\mathcal{R}(k)$ on a $16{\tilde{M}}^3$ grid. The reason becomes clear if we consider the formulae in section 4.3 in \cite{fsewald2017}. It is apparent that we  need to truncate the  Green's function values centered at a particular point in the physical grid. Thus, to perform the truncation, we start with sampling values of $\widehat H^\mathcal{R}(k)$, perform an IFFT to obtain values in physical space, and then perform an FFT after truncation.

\begin{algorithm}[H]
\caption{SE method for Fourier-space calculation of stokeslet for \halfspace}
\label{alg:SE_halfspace}
\begin{algorithmic}[1]
\STATE \textbf{Input}: Source locations $\mb{y}^{(m)} \in [0, L)^3$, source intensities $\mb{f}^{(m)}, \;$ for $m=1,2, \dotsc N$. 
Ewald parameter $\xi$, Gaussian support width $P$, basic grid size $M$ so that grid spacing $h= L/M$.
\STATE {Extend domain to $[0,L) \times [0, L) \times [-L, L)$ and  set wall at $x_3 = 0$.} Now construct reflections ${\mb{y}^{(m)}}^I, \;  {\mb{f}^{(m)}}^I$ for $m=1 \dotsc N$. 
\STATE Extend domain to adjust for Gaussian support $\tilde{L} = L + Ph \Rightarrow \tilde{M} = M +P$. Thus uniform grid of dimension $\tilde{M} \times \tilde{M} \times 2\tilde{M}$ on domain $\tilde{L} \times \tilde{L} \times 2\tilde{L}$.
\STATE Set $\eta = (P\xi^2h^2)/(c^2\pi) > 0$ where $c=0.95$.
\STATE Evaluate the sum of truncated Gaussians $C_1$ for stokeslet and image and $C_2, C_3$  where 
\begin{align*}
& C_1(\mb{x})= \sum_{m=1}^{2N} \left(\dfrac{2\xi^2}{\pi\eta}\right)^{3/2}
  \left[ ( \mb{F}^{(m)}  ) e^{-2\xi^2\norm{ \mb{x} - \mb{z}^{(m)} }/ \eta}  \right], \\
 & C_2(\mb{x})= \sum_{m=1}^N \left(\dfrac{2\xi^2}{\pi\eta}\right)^{3/2} \left( f^{(m)}_3 \right) e^{-2\xi^2\norm{ \mb{x} - {\mb{y}^{(m)}}^I }/ \eta}, \\
 & C_3(\mb{x})= \sum_{m=1}^N \left(\dfrac{2\xi^2}{\pi\eta}\right)^{3/2} \left( y^{(m)}_3 {\mb{f}^{(m)} }^I \right) e^{-2\xi^2\norm{ \mb{x} - {\mb{y}^{(m)}}^I }/ \eta}.
\end{align*}
Here $\mb{F}^{(m)}, f^{(m)}_3, y^{(m)}_3 {\mb{f}^{(m)} }^I$ play the role of {$\mb{c}^{(m)}$} in the general formula \eqref{eqn:uF}.
\STATE Compute the FFTs, i.e., $\widehat{C}_1(\mb{k}), \widehat{C}_2(\mb{k}), \widehat{C}_3(\mb{k})$, zero-padded by a factor of 2. 
\STATE Use precomputed $\widehat{H}^{\mathcal{R}}(k), \widehat{B}^{\mathcal{R}}(k)$ with $\mathcal{R} = \sqrt{6}\tilde{L}$ oversampled by factor of 2 (see remarks on precomputation). This corresponds to {$\widehat{\mathcal{K}}(k)$} in \eqref{eqn:uF}. 
\STATE Compute the product $\mb{A}(\mb{k}, \xi) \widehat{\mathcal{K}}(k) e^{-(1-\eta)\frac{k^2}{4\xi^2}}  {\widehat{C}(\mb{k})}$ with reference to the general representation \eqref{eqn:uF} appropriate  for each of the 5 integrals, and where $\widehat{C}(\mb{k})$ hence is one of $\widehat{C}_i(\mb{k}), i=1, 2, 3$. Call them $\widehat{\mb{t}_1}, \widehat{\mb{t}_2}, \widehat{\mb{t}_3}, \widehat{t_4},\widehat{t_5}$ The first three are vectors and the last two are scalars.
\STATE Compute the IFFT of $\widehat{\mb{t}_1} + (0, 0, \widehat{t_4}+\widehat{t_5})^T$ and 
$\widehat{\mb{t}_2} + \widehat{\mb{t}_3}$, truncate to get $\mb{T}_1, \; \mb{T}_2$ respectively on $\tilde{M} \times \tilde{M} \times 2\tilde{M}$ grid.
\STATE By using trapezoidal rule and truncated Gaussians at the target point {$\mb{x} = \mb{y}^{(m)}$}, compute {
\begin{align*}
\mb{u}_F(\mb{y}^{(m)}, \xi) =& \frac{1}{(2\pi)^3} \int_{\mathbb R^3} \mb{T}_1(\mb{y})\left( \frac{2\xi^2}{\pi\eta}\right)^{3/2} e^{-2\xi^2\norm{ \mb{y}^{(m)} - \mb{y} }/ \eta} d\mb{y} \\
 & - x_3 \frac{1}{(2\pi)^3} \int_{\mathbb R^3}\mb{T}_2(\mb{y})  \left( \dfrac{2\xi^2}{\pi\eta}\right)^{3/2} e^{-2\xi^2\norm{ \mb{y}^{(m)} - \mb{y} }/ \eta} d\mb{y}.
 \end{align*}
}
\STATE \textbf{Output}: $\mb{u}_F(\mb{y}^{(m)}, \xi), \; m = 1 \dotsc N$. 
\end{algorithmic}
\end{algorithm}

The computation of $\mb{u}_F(\mb{x}, \xi)$ has been outlined in  Algorithm~\ref{alg:SE_halfspace} and can be organized into the following main steps:
\begin{enumerate}
\item \emph{Preprocessing} (Steps 2, 3): Setting up the appropriate computational domain, performing reflections etc.
\item \emph{Spreading} (Step 5): Computing the data on the grid using truncated Gaussians. This is  essentially the source term or some component(s) of it scaled in various ways. The Gaussians are assumed to have a support of $P^3$ grid points. {(Explicit formula given in Algorithm~\ref{alg:SE_halfspace})}
\item \emph{FFT} (Step 6): Computing the three-dimensional FFT, zero-padded to double the size. The factor of 2 is necessary to perform the aperiodic convolution. Note that for the stokeslet, we perform 2 vector FFTs and one scalar FFT.
\item \emph{Precomputation} (Step 7) Computing $\widehat H^\mathcal{R}(k)$ or $\widehat B^\mathcal{R}(k)$. This step depends on the size of the computational domain only and is actually performed after Steps 2 and 3, but is listed later here for aiding the flow of ideas.
\item \emph{Scaling} (Steps 8 and 9): Computing product of FFT with  $\mb{A}(\mb{k}, \xi)$ etc. through precomputed $\widehat H^\mathcal{R}(k)$ or $\widehat B^\mathcal{R}(k)$. In this step all the quantities involved have a twice oversampled representation.
\item \emph{IFFT} (Step 10): Applying the inverse three-dimensional FFT to the result of the scaling and truncating to obtain result on  $\tilde{M} \times \tilde{M} \times 2\tilde{M}$ grid.
\item \emph{Quadrature} (Step 11): Evaluating the resulting integral by trapezoidal rule for each $\mb{x}$.
\end{enumerate}

In comparison to the evaluation for \freespace in \cite{fsewald2017}, the number of FFTs and IFFTs increases and this information is summarized in Table~\ref{table:fse-hse}.

\begin{table}[ht]
\centering
\caption{Breakdown of the Fourier-space operations. Note that the  kernel-image sums are combined for the purpose of FFT into a single consolidated sum over $2N$ terms  by modifying the image-sources.}
\label{table:fse-hse}
\begin{tabular}{|l|c|c|c|c|c|c|}
\hline
& \multicolumn{2}{c}{\textbf{Stokeslet}}   & \multicolumn{2}{|c|}{\textbf{Stresslet}}      
   & \multicolumn{2}{c|}{\textbf{Rotlet}}   \\ \hline
& \multicolumn{1}{|l|}{FFT} & \multicolumn{1}{l}{IFFT} &  
\multicolumn{1}{|l|}{FFT} & \multicolumn{1}{l}{IFFT} & 
\multicolumn{1}{|l|}{FFT} & \multicolumn{1}{l|}{IFFT}  \\ \hline
\textbf{Free-space} & 3     & 3    & 9    & 3    & 3 & 3  \\
\textbf{Half-space} & 7     & 6    & 21   & 6    & 6 & 6   \\ \hline
\end{tabular}
\end{table}

The increase in the number of FFTs is consistent with the \halfspace formula 
\eqref{eqn:v_sum_halfspace}; the correction terms require one vector and one scalar FFT (3 + 1) for the stokeslet, one vector FFT (3) for the rotlet, and  one tensor and one vector FFT (9 + 3) for the stresslet.
The increase in the number of FFTs required is thus not surprising, but the increase in the number of IFFTs is not so obvious, for after calculating the FFTs and performing the convolution as a product in Fourier space, it might seem that we could combine them all and perform a single IFFT and integration step. However, that is not possible due to the presence of the coordinate $x_3$ multiplying the gradient. This  term needs to be treated separately and hence it requires an extra IFFT and integration step. 
One might well ask if the Ewald summation technique is still worthwhile with this increase in the complexity of the problem in Fourier space; fortunately, the answer is affirmative, as demonstrated in the next section.

\section{Truncation Errors}
\label{sec:truncation}

The errors in the real-space calculation are caused by truncation. However, the errors in the SE method for the Fourier-space part are not due to truncation alone; there are approximation errors due to the quadrature rule for integration and the discretization and truncation of the Gaussians in the spreading and quadrature steps. Given $P$, the number of points across a truncated Gaussian, the parameter $\eta = \eta(P)$ can be chosen to balance discretization and truncation errors, which leaves $P$ as the only free parameter. Approximation errors decay exponentially with $P$ and by using a sufficiently large $P$, the approximation errors may be considered negligible so that the measured error is due to truncation errors only \cite{Lindbo2011c}, as the Ewald parameter $\xi$ does not introduce any errors.

The truncation errors appear in the real-space part because we consider only local interactions within the radius $r_c$ while in the Fourier-space part, the integral in \eqref{eqn:uF} over all of $\mathbb{R}^3$ is truncated in practice to  consider some large but finite wave number with magnitude $k_{\infty}$. The value of $k_{\infty}$ is related to the grid and computational domain by the relation $k_{\infty} = \pi/h = \pi\tilde{M}/\tilde{L}$.
The \emph{exact} real-space contribution is obtained by letting $r_c \rightarrow \infty$ in the real-space sum, and in combination with the naive direct sum, it also yields the \emph{exact} Fourier-space contribution\footnote{by \emph{exact} Fourier-space contribution, we mean the integral in \eqref{eqn:uF} evaluated over all of $\mathbb{R}^3$}. From these, the computed truncation errors are found. In \cite{fsewald2017}, the authors report truncation error estimates based on the methodology introduced by Kolafa \& Perram \cite{Kolafa1992} that agree closely with the computed errors for both real-space and Fourier-space evaluations. These are statistical error estimates for the root mean square (RMS) truncation error, defined as
  $$ \delta \mb{u} = \sqrt{ \frac{1}{N} \sum_{m=1}^N \norm{ \mb{
        u}_{\mathrm{exact}}(\mb{y}^{(m)}) - \mb{u}_{\mathrm{computed}}(\mb{y}^{(m)}) }^2 } .$$
 Implicit in these estimates are the assumptions that the sources are randomly distributed, and that the error measure has a Gaussian distribution. 

Since the derived estimates are statistical in nature, the methodology of Kolafa \& Perram will not be able to account for possible cancellations due to symmetry in the \halfspace formula.  On the contrary, this approach will yield an estimate that is the sum of two contributions--
\begin{enumerate}
\item Truncation errors due to the stokeslet, stresslet or rotlet kernels (as in the Tables)
\item Truncation errors due to the correction terms 
\end{enumerate}
Such an estimate is likely to be a conservative upper-bound for the computed errors. Therefore, we examine the truncation error expressions to ascertain if one can neglect the error contribution made by the correction terms. 

From Tables~\ref{tab:est_fourier} \& \ref{tab:est_real}, it is evident that the kernel estimates have an exponential decay term ($e^{-\xi^2 r_c^2}$ or $e^{-\kmax^2/4\xi^2}$) multiplied by some power of $r_c$ or $k_{\infty}$. In case of the \halfspace, the correction terms are all harmonic functions and/or their derivatives,  and their truncation error contribution (see \cite{Kolafa1992}) decays much faster than the stokeslet or stresslet kernels and at the same rate as the rotlet kernel. Thus, the existing truncation error estimates for \freespace are a good starting point for the \halfspace as well. Of course, in evaluating the estimate, the sum is over both {sources and images}, and this modifies the quantity $Q$ and the RMS error.

\begin{figure}[htbp]
  \centering
  \begin{subfigure}[b]{0.32\textwidth}
    \centering
    \includegraphics[width=\textwidth]{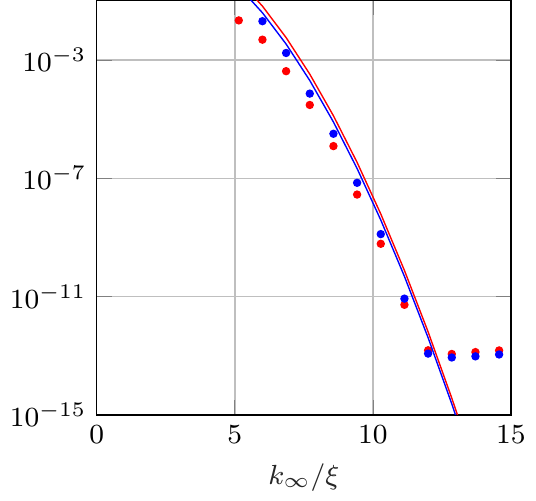}
    \caption{Stokeslet}
    \label{fig:fourier_err:stokeslet}
  \end{subfigure}  
  \begin{subfigure}[b]{0.32\textwidth}
    \centering
    \includegraphics[width=\textwidth]{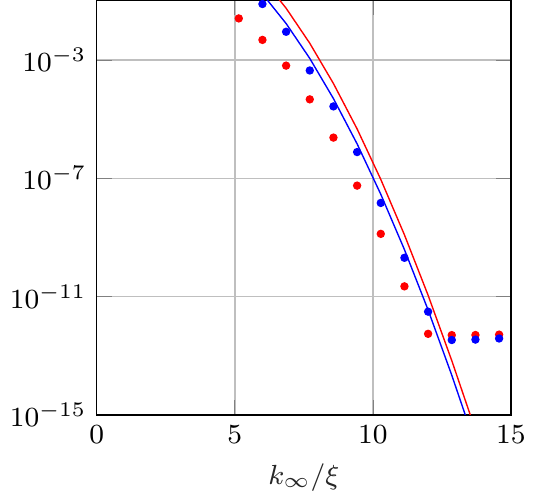}
    \caption{Stresslet}
    \label{fig:fourier_err:stresslet}
  \end{subfigure}  
  \begin{subfigure}[b]{0.32\textwidth}
    \centering
    \includegraphics[width=\textwidth]{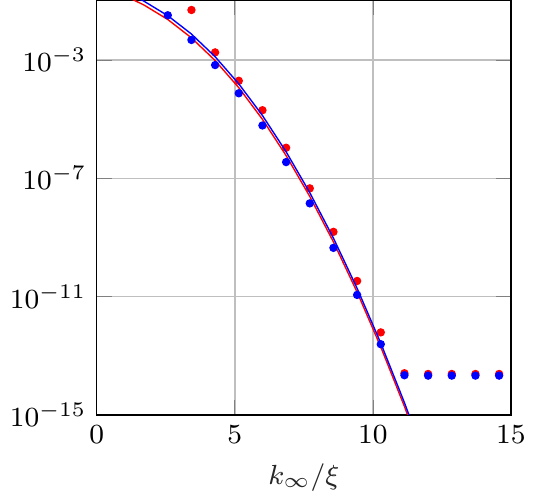}
    \caption{Rotlet}
    \label{fig:fourier_err:rotlet}
  \end{subfigure}  
  \caption{RMS of relative Fourier-space truncation errors for the
    stokeslet, stresslet and rotlet. Dots are measured value, solid
    lines are estimates based on Table~\ref{tab:est_fourier}. Red coloured dots and lines are associated with the \halfspace while blue colour is associated with the \freespace  evaluation. The system is $N=10^4$ randomly distributed
    point sources in a cube with sides $L=3$, with $k_{\infty}=\pi \tilde{M}
    /\tilde{L}$, $\xi=3.49$, $M=1\dots50$.}
  \label{fig:fourier_err}
\end{figure}

\begin{figure}[htbp]
  \centering
  \begin{subfigure}[b]{0.32\textwidth}
    \centering
    \includegraphics[width=\textwidth]{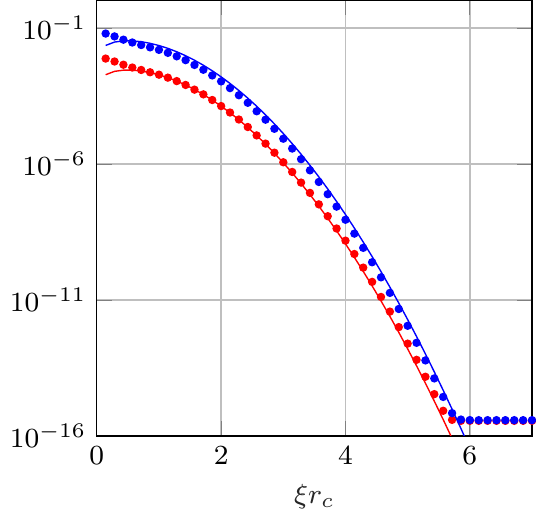}    
    \caption{Stokeslet}
    \label{fig:real_err:stokeslet}
  \end{subfigure}  
  \begin{subfigure}[b]{0.32\textwidth}
    \centering
    \includegraphics[width=\textwidth]{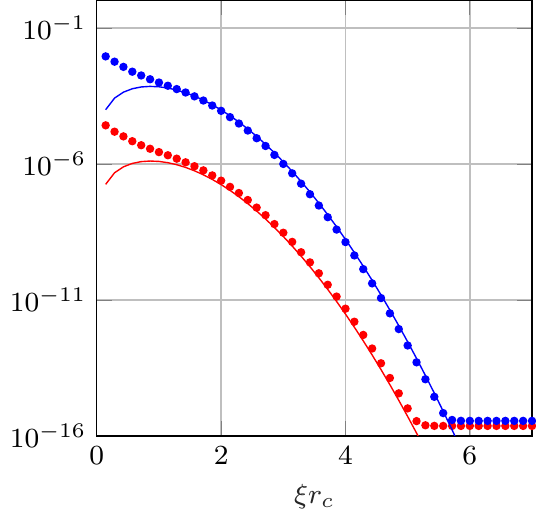}
    \caption{Stresslet}
    \label{fig:real_err:stresslet}
  \end{subfigure}  
  \begin{subfigure}[b]{0.32\textwidth}
    \centering
    \includegraphics[width=\textwidth]{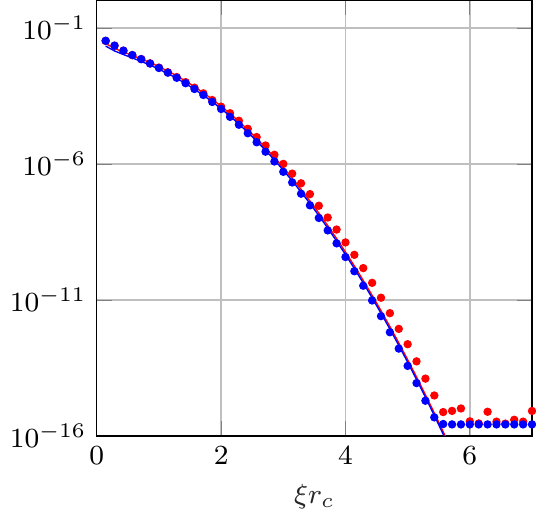}    
    \caption{Rotlet}
    \label{fig:real_err:rotlet}
  \end{subfigure}  
  \caption{RMS of relative real-space truncation errors for the
    stokeslet, stresslet and rotlet. Dots are measured value, solid
    lines are computed using the estimates of Table~\ref{tab:est_real}. Red coloured dots and lines are associated with the \halfspace while blue colour is associated with the 
    \freespace  evaluation.The system is $N=2000$ randomly distributed
    point sources in a cube with sides $L=3$, with $\xi=4.67$, and
    $r_c \in [0, L/2]$.}
  \label{fig:real_err}
\end{figure}

In Figures~\ref{fig:fourier_err} and \ref{fig:real_err}, the truncation error estimate is compared to the computed error for both real-space and Fourier-space parts. This is done by calculating the relative RMS error which is the ratio of the RMS error and the RMS value of the velocity considered at all targets.   For reference, the corresponding  curves for the \freespace problem are also plotted. For both the stokeslet and stresslet in Figures~\ref{fig:fourier_err} and \ref{fig:real_err}, the computed errors are lesser than the \halfspace estimate, and  the  \freespace computed error and estimate. The computed error for the \halfspace here is lower due to the cancellations induced by symmetry. For the rotlet, while all these curves are much closer to each other, closer examination reveals that the previous trend is no longer upheld. This is because the correction terms that have been neglected have the same order as the kernel. However, the overall agreement of the computed error with the existing estimates justifies the decision to neglect the contribution of the correction terms to the truncation error estimate.

\begin{table}[htbp]
  \centering
  \begin{tabular}{c|c|c|c}
    & Stokeslet, $\stokeslet^F$ 
    & Stresslet, $\stresslet^F$ 
    & Rotlet, $\rotlet^F$ 
    \\ \hline
    &&&\\
    $\delta\mb{u}^F$ & 
    $\displaystyle \sqrt{Q}\frac{\trunc \kmax^3}{\xi^2 \pi L} e^{-\kmax^2/4\xi^2}$  &
    $\displaystyle \sqrt{\frac{7Q}{6}}\frac{\trunc \kmax^4}{\xi^2 \pi L} e^{-\kmax^2/4\xi^2}$ &
    $\displaystyle \sqrt{\frac{8\xi^2Q}{3\pi L^3 \kmax}} e^{-\kmax^2/4\xi^2}$ 
    \\
    $Q$ & $\displaystyle \sum_{m=1}^N \norm{ \mb{f}^{(m)} }^2$ & $\displaystyle \sum_{m=1}^N \norm{ \mb{g}^{(m)} \otimes \mb{q}^{(m)} }^2_{F}$ & $\displaystyle \sum_{m=1}^N \norm{ \mb{g}^{(m)} \times \mb{q}^{(m)} }^2$ \\
    &&&\\
  \end{tabular}
  \caption{Fourier-space truncation errors for the stokeslet, stresslet, and rotlet \cite{AfKlinteberg2016rot} for \freespace. The quantity $Q$ is defined for each kernel as in the second row and $\norm{\cdot}_F$ denotes the Frobenius norm for second-order tensors.}
  \label{tab:est_fourier}
\end{table}

\begin{table}[htbp]
  \centering
  \begin{tabular}{c|c|c|c}
    & Stokeslet, $\stokeslet^R$ 
    & Stresslet, $\stresslet^R$ 
    & Rotlet, $\rotlet^R$ 
\\
    \hline
    &&&\\
    $\delta\mb{u}^R$ 
    & $\displaystyle \sqrt{\frac{4 Q r_c}{L^3}} e^{-\xi^2 r_c^2}$
    & $\displaystyle \sqrt{\frac{112 Q \xi^4 r_c^3}{ 9 L^3 }} e^{-\xi^2 r_c^2}$ 
    & $\displaystyle \sqrt{\frac{8 Q}{3 L^3 r_c}} e^{-\xi^2 r_c^2}$ 
    \\
    &&&\\
  \end{tabular}
  \caption{Real-space truncation errors for the stokeslet \cite{Lindbo2011e}, stresslet, and rotlet \cite{AfKlinteberg2016rot} for \freespace. The quantity $Q$ is defined in the same way as for the Fourier component in the previous table.}
  \label{tab:est_real}
\end{table}

\section{Numerical Results}
\label{sec:results}

We consider $N$ random point sources drawn from a uniform distribution from a box  of dimension $L \times L \times L$. The sum \eqref{eqn:v_sum_free_space} is evaluated with stokeslets, stresslets and rotlets, at the same $N$ target locations.  All components of the
source strengths and source orientations are random numbers from a uniform distribution
on $[-1,1]$.  All computationally intensive routines are written in C
and are called from Matlab using MEX interfaces.  The results are obtained on a laptop workstation with an Intel Core i7-6500U Processor (2.50 GHz) and 16 GB of memory, running two cores unless stated so.  Actual errors are measured by comparing the
result with evaluating the sum by direct summation.

For a given system of $N$ charges in a box of edge length $L$, the required
parameters for our \halfspace Ewald method are the Ewald parameter
$\xi$, the real space truncation radius $r_c$, the number of grid
points $M$ covering the computational domain, and the Gaussian support
width $P$. Other parameters like $\delta_L, \tilde{L}, \tilde{M}$ are then set automatically from these (see Algorithm~\ref{alg:SE_halfspace} and \cite{fsewald2017} for details). For a large-scale numerical computation, $\xi$, $r_c$, $M$ and $P$ must be set optimally. For a given value of $\xi$ and absolute error tolerance $\epsilon$ for the \freespace evaluation, close-to optimal values for $M$ and $r_c$ were computed in \cite{fsewald2017}  using the truncation error estimates in Tables ~\ref{tab:est_fourier} and
\ref{tab:est_real}. We use the very same optimal values as starting points for the \halfspace since the same truncation error estimates hold good for the \halfspace too as shown. Then we  perturb $\xi$ to  achieve the smallest runtime while keeping $\xi r_c$ and $\frac{k_{\infty}}{\xi}$ constant. Note that $k_{\infty} = \frac{\pi \tilde{M}}{\tilde{L}}$. The optimal set of values found is used for larger systems by scaling $k_{\infty}$ such that $\frac{L}{k_{\infty}}$ is constant.

\begin{figure}[htb]
\centering
  \includegraphics[width=0.6\textwidth, height=0.4\textwidth]{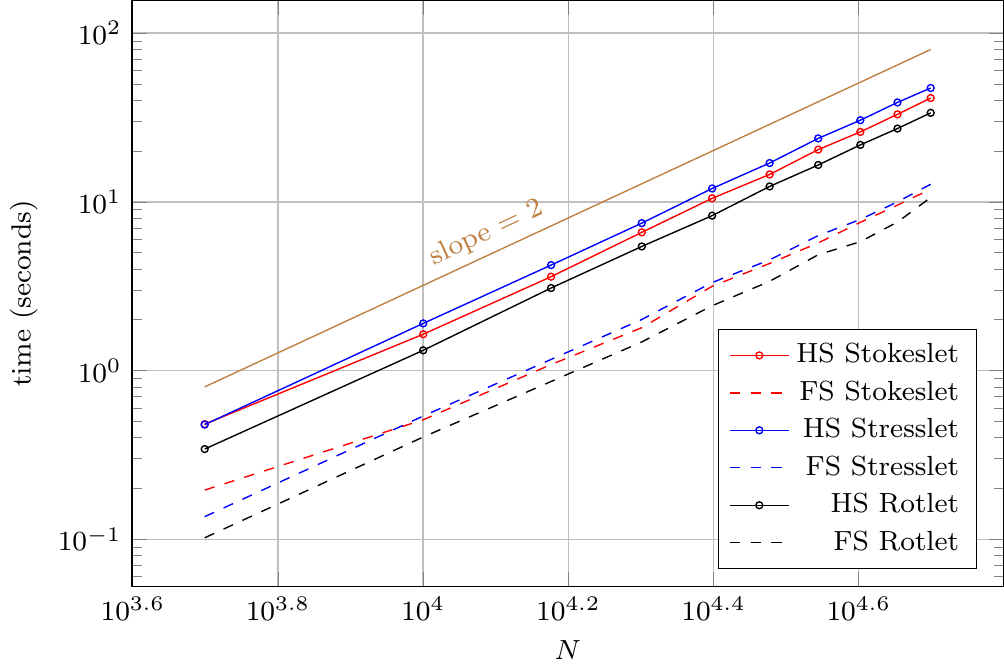}
  \caption{CPU runtimes of the direct sum for the \halfspace (HS) and \freespace (FS) as a function of $N$. Both are $O(N^2)$, and  the spacing indicates their ratio, which is in the range 3.3--3.7.}
  \label{fig:direct_sum_HS_FS}  
\end{figure}

\begin{figure}[htb]
\centering
\begin{subfigure}[b]{0.48\textwidth}
    \centering
    \includegraphics[width=\textwidth, height= 0.8\textwidth]{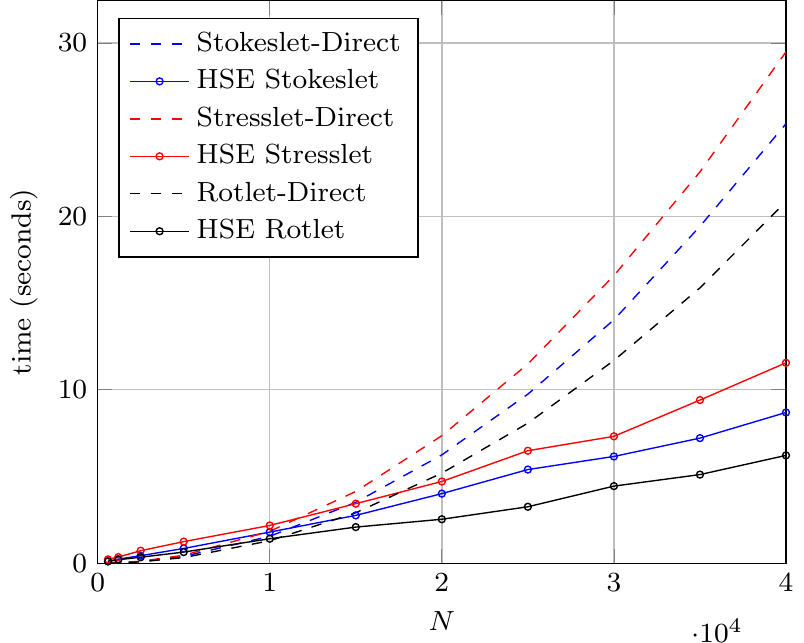}    
    \caption{Half-Space}
    \label{fig:CPU_runtime_HS}
 \end{subfigure} 
 \begin{subfigure}[b]{0.48\textwidth}
    \centering
    \includegraphics[width=\textwidth, height=0.8\textwidth]{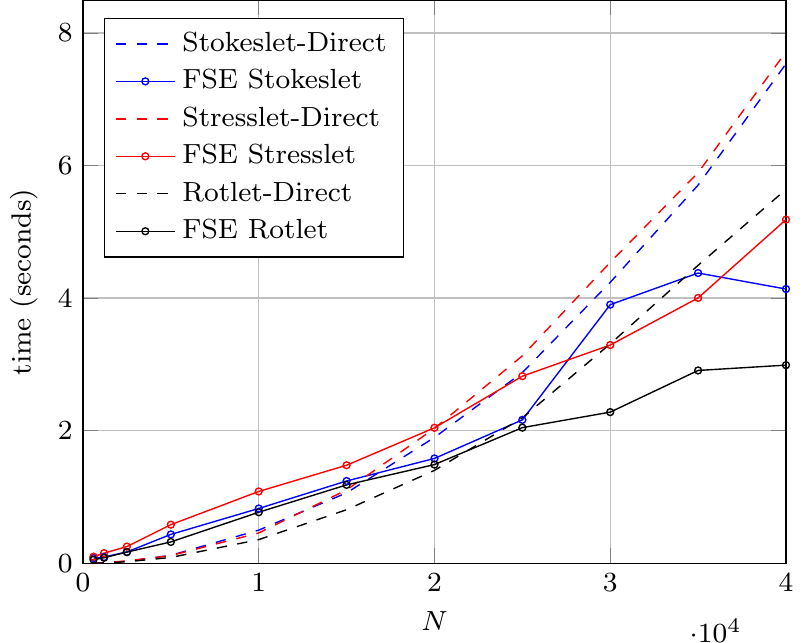}    
    \caption{Free-Space}
    \label{fig:CPU_runtime_FS}
 \end{subfigure} 
\caption{CPU run-time for direct and Ewald sum with different kernels}
\label{fig:CPU_runtime}
\end{figure}

The results for \freespace in \cite{fsewald2017} convincingly demonstrated the need for the Ewald decomposition by exhibiting  the speed-up gained  over the  naive direct sum. The  present work aims to make a similar case for the \halfspace. Before that however, it would be interesting to compare the computational expense of evaluating the direct sums themselves for the \freespace and the \halfspace formulae. This has been done in Figure~\ref{fig:direct_sum_HS_FS} and it is seen that on average, the \halfspace sum is between 3.3--3.7  times more expensive (with rotlet the least, and stresslet the most) than the \freespace sum and this  factor stays constant even as the number of sources increases.  The obtained range 3.3--3.7 is reasonable since the direct sum for the \halfspace  involves $4N$ terms (kernel and correction terms), in contrast to that for the \freespace which involves only $N$ terms. The additional expense due to the correction terms is smallest for the rotlet, and most for the stresslet, and this is not surprising seeing the formulae \eqref{eqn:phi}.
We now tackle the question of whether it is worthwhile to consider the Spectral Ewald method  for the \halfspace formula despite the substantial increase in the number of FFTs and IFFTs that need to be performed for the Fourier-space component. In Figure~\ref{fig:CPU_runtime_HS}
the computing time for evaluation of the sums is plotted versus $N$, for all three kernels and for both the Half Space Ewald (HSE) method and direct summation.  The relative RMS error
for HSE is kept below $0.5 \times 10^{-8}$.  As we vary $N$, we maintain a constant density
$N/L^3=2500$ by changing the size of the box. The system is thus scaled while other parameters of the method like  $\xi, r_c$, $P$ and the grid resolution $L/M$ are kept constant. As  both $N$ and $L$ increase, so does the grid size.  

For all kernels, we set $P = 16$, while the pair $(\xi, r_c)$ is (6, 0.76), (5.8, 0.76), (7, 0.63) for the stokeslet, stresslet, and rotlet respectively. When $L=2$, $M=40$, $36$ and $38$ for the three kernels, and the ratio $M/L$ is kept constant as the system is scaled.
We have excluded the precomputation cost in our  runtimes since it is performed only once, and is easily amortized over multiple runs due to iterations or time steps when the size of the domain does not change. The figure allows us to determine the break-even value $N$, that is, the smallest value of $N$ for which the Ewald summation is faster than the direct sum. 
In order to compare CPU runtimes, it is necessary that the simulations should use the same system and number of cores. Hence the CPU runtime study for \freespace was repeated with the same parameters as above except that $\xi=7$, and $r_c = 0.63, 0.63, 0.58$ for the stokeslet, stresslet and rotlet respectively. The results are shown in Figure~\ref{fig:CPU_runtime_FS}.

The break-even values for all kernels (when excluding the precomputation cost) obtained from Figure~\ref{fig:CPU_runtime} are presented in Table~\ref{tab:break_even}. Since the cost of the direct sum for the \halfspace almost quadruples in comparison to that of \freespace, it benefits more from the Ewald decomposition and Spectral Ewald method and the break-even is attained much earlier. As expected, the stresslet, with the steep increase in the number of FFTs, has the largest break-even, but it is only slightly greater  than the  break-even value for the other kernels. These numbers underline the advantage of the Spectral Ewald method for the evaluation of the formulae \eqref{eqn:v_sum_halfspace} for the \halfspace.
\begin{table}[htb]
\centering
\begin{tabular}{|c|c|c|}
\hline
\emph{Kernel}   & \emph{Half Space Ewald}    & \emph{Free Space Ewald} \\ \hline      Stokeslet & $1.15 \times 10^{4}$ & $1.7 \times 10^{4}$ \\ 
Stresslet & $1.2 \times 10^{4}$ & $2.00 \times 10^{4}$  \\ 
Rotlet    & $1.05 \times 10^{4}$ & $2.2 \times 10^{4}$   \\ \hline               
\end{tabular}
\caption{Number of sources $N$ at which break-even run times occur with Ewald method, with the relative RMS error under $0.5 \times 10^{-8}$}
\label{tab:break_even}
\end{table}

\begin{figure}[htb]
\centering
\includegraphics[width=0.49\textwidth, height=0.4\textwidth]{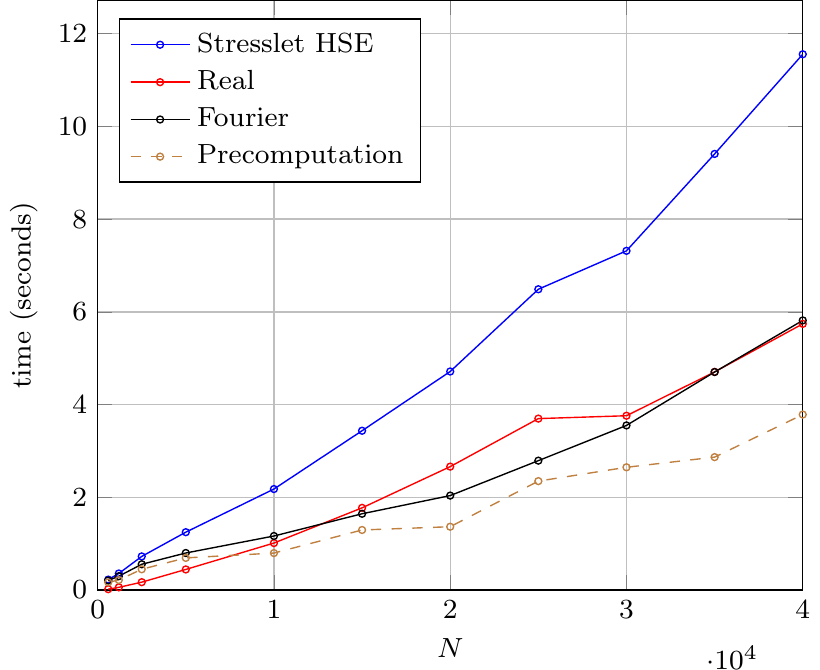}
\includegraphics[width=0.49\textwidth, height=0.4\textwidth]{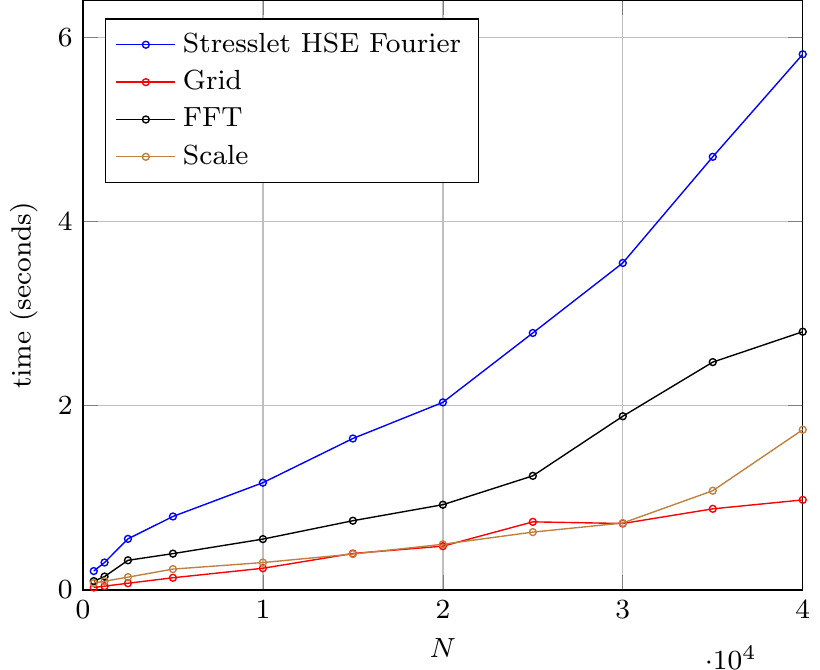}
\caption{Breakdown of runtimes (left) and Fourier-space runtime
    (right) for evaluating the stresslet as a function of number of
    particles. The  runtime for precomputation is also shown in the left plot, this is excluded from the overall runtime. }
\label{fig:CPU_run_time_breakdown}
\end{figure}

We next study the computational run-time of different parts of the algorithm. In the left  plot of Figure~\ref{fig:CPU_run_time_breakdown} the runtimes for real-space and Fourier-space evaluation are presented for the stresslet\footnote{We choose the stresslet because it is the most complicated}, along with the cost of precomputation of the Green's functions. For the stresslet kernels, the precomputation involves the  evaluation of the Green's function for the stresslet as well as that of the correction functions and including it in the total runtime will increase the cost  by 25--33\%.
The usage of an  optimized value of the Ewald parameter $\xi$ balances the cost of the real-space and Fourier-space evaluation and they are thus of the same order of magnitude. The right plot of Figure~\ref{fig:CPU_run_time_breakdown} shows a further breakdown of the Fourier-space cost into its main  constituent steps, namely, \emph{Gridding}, \emph{Scaling}, and \emph{FFT}. The  scaling step is clearly the cheapest among them,  and the overall
results are very similar to the case of  \freespace despite the fact that here we perform 21 FFTs compared to only 9 for \freespace.

\section{Conclusion and further work}
We have presented a fast summation method for the \halfspace Green's functions of Stokes flow derived by Gimbutas \etal \cite{Gimbutas2015}. The fast summation method and its implementation follows that  for \freespace Green's functions by Klinteberg \etal \cite{fsewald2017}. The method is based on the Ewald decomposition that recasts the sum into a sum of two exponentially decaying series: one in real-space (short-range interactions) and one in Fourier-space (long-range interactions) with the convergence of each series  controlled by a common parameter.

While the evaluation of the real-space component proceeded along expected lines, the presence of extra terms complicated the task for the Fourier-space component. We followed the framework of the Spectral Ewald method for \freespace Stokes flow introduced recently, and exploited the structure of the terms to optimize the number of FFTs and IFFTs that need to be performed. Furthermore, we demonstrated that with very elementary modifications the truncation error estimates for \freespace Stokes flow remain valid. 

The implementation for the \halfspace does incur extra costs in comparison to the \freespace in  multiple ways such as the gridding of a larger computational domain, substantial increase in the number of FFTs and IFFTs that need to be evaluated but the computational savings are also greater.

Future work can take shape in one of two ways. For one, it would be beneficial to use this work in the framework of a boundary integral method for Stokes flow in a \halfspace motivated by a physical problem of sedimentation. 
A more involved, and mathematically interesting question would be to consider a 2-periodic extension with periodicity in the in-plane directions, akin to the 2-periodic extension for Stokes flow considered by Lindbo \etal \cite{Lindbo2011e}.

\section{Acknowledgements}
  Anna-Karin Tornberg thanks the G\"{o}ran Gustafsson Foundation
  for Research in Natural Sciences and Medicine and the Swedish
    e-Science Research Centre (SeRC). Shriram Srinivasan gratefully
  acknowledges the financial support of Linn{\'e} Flow Centre and thanks Ludvig af Klinteberg for helpful ideas and advice with the numerical implementation.


\bibliography{library}

\begin{thebibliography}{10}

\bibitem{Gimbutas2015}
Z.~Gimbutas, L.~Greengard, and S.~Veerapaneni.
\newblock Simple and efficient representations for the fundamental solutions of
  {S}tokes flow in a half-space.
\newblock {\em Journal of Fluid Mechanics}, 776, 2015.

\bibitem{fsewald2017}
Ludvig af~Klinteberg, Davoud-Saffar Shamshirgar, and Anna-Karin Tornberg.
\newblock Fast {E}wald summation for free-space {S}tokes potentials.
\newblock {\em Research in the Mathematical Sciences}, 4(1), 2017.

\bibitem{KimKarilla1991}
Sangtae Kim and Seppo~J. Karilla.
\newblock {\em Microhydrodynamics}.
\newblock Butterworth-Heineman, 1991.

\bibitem{JBlake1971}
J.~R. Blake.
\newblock A note on the image system for a stokeslet in a no-slip boundary.
\newblock {\em Mathematical Proceedings of the Cambridge Philosophical
  Society}, 1971.

\bibitem{Ewald1921}
P.~P. Ewald.
\newblock {Die Berechnung optischer und elektrostatischer Gitterpotentiale}.
\newblock {\em Annals of Physics}, 369(3):253--287, 1921.

\bibitem{Hasimoto1959}
H.~Hasimoto.
\newblock {On the periodic fundamental solutions of the Stokes equations and
  their application to viscous flow past a cubic array of spheres}.
\newblock {\em Journal of Fluid Mechanics}, 5(02):317, feb 1959.

\bibitem{Deserno1998}
Markus Deserno and Christian Holm.
\newblock {How to mesh up Ewald sums. I. A theoretical and numerical comparison
  of various particle mesh routines}.
\newblock {\em Journal of Chemical Physics}, 109(18):7678, 1998.

\bibitem{Essmann1995}
U.~Essmann, L.~Perera, M.~L. Berkowitz, T.~Darden, H.~Lee, and L.~G. Pedersen.
\newblock {A smooth particle mesh Ewald method}.
\newblock {\em Journal of Chemical Physics}, 103(19):8577--8593, 1995.

\bibitem{Saintillan2005}
D.~Saintillan, E.~Darve, and E.~Shaqfeh.
\newblock {A smooth particle-mesh Ewald algorithm for Stokes suspension
  simulations: The sedimentation of fibers}.
\newblock {\em Physics of Fluids}, 17(3), 2005.

\bibitem{Lindbo2010}
Dag Lindbo and Anna-Karin Tornberg.
\newblock {Spectrally accurate fast summation for periodic Stokes potentials}.
\newblock {\em Journal of Computational Physics}, 229(23):8994--9010, 2010.

\bibitem{AfKlinteberg2014a}
Ludvig {Af Klinteberg} and Anna-Karin Tornberg.
\newblock {Fast Ewald summation for Stokesian particle suspensions}.
\newblock {\em International Journal for Numerical Methods in Fluids},
  76(10):669--698, dec 2014.

\bibitem{AfKlinteberg2016rot}
Ludvig af~Klinteberg.
\newblock {Ewald summation for the rotlet singularity of Stokes flow}.
\newblock {\em arXiv:1603.07467 [physics.flu-dyn]}, mar 2016.

\bibitem{Lindbo2011e}
Dag Lindbo and Anna-Karin Tornberg.
\newblock {Fast and spectrally accurate summation of 2-periodic Stokes
  potentials}.
\newblock {\em arXiv:1111.1815v1 [physics.flu-dyn]}, 2011.

\bibitem{SE1P_Ewald_electrostatics}
Davoud-Saffar Shamshirgar and Anna-Karin Tornberg.
\newblock {The Spectral Ewald method for singly periodic domains}.
\newblock {\em Journal of Computational Physics}, 347:341--366, 2017.

\bibitem{Vico2016}
Felipe Vico, Leslie Greengard, and Miguel Ferrando.
\newblock {Fast convolution with free-space Green's functions}.
\newblock {\em Journal of Computational Physics}, 323:191--203, oct 2016.

\bibitem{Greengard1987}
L.~Greengard and V.~Rokhlin.
\newblock {A Fast Algorithm for Particle Simulations}.
\newblock {\em Journal of Computational Physics}, 73:325--348, 1987.

\bibitem{Rodin2000}
Y.~Fu and G.~J. Rodin.
\newblock {Fast solution methods for three-dimensional Stokesian many-particle
  problems}.
\newblock {\em Communications in Numerical Methods in Engineering},
  16:145--149, 2000.

\bibitem{Duraiswami2006}
N.~A. Gumerov and R.~Duraiswami.
\newblock {Fast multipole method for the biharmonic equation in three
  dimensions}.
\newblock {\em Journal of Computational Physics}, 215:363--383, 2006.

\bibitem{Tornberg2008}
Anna-Karin Tornberg and Leslie Greengard.
\newblock {A fast multipole method for the three-dimensional Stokes equations}.
\newblock {\em Journal of Computational Physics}, 227(3):1613--1619, jan 2008.

\bibitem{papkovich1932solution}
P.~F. Papkovich.
\newblock Solution g{\'e}n{\'e}rale des {\'e}quations differentielles
  fondamentales d’{\'e}lasticit{\'e} exprim{\'e}e par trois fonctions
  harmoniques.
\newblock {\em Comptus Rendus de l' Acad{\'e}mie de Sciences. Paris},
  195(3):513--515, 1932.

\bibitem{neuber1934neuer}
H.~Neuber.
\newblock Ein neuer ansatz zur l{\"o}sung r{\"a}umlicher probleme der
  elastizit{\"a}tstheorie. der hohlkegel unter einzellast als beispiel.
\newblock {\em Zeitschrift f{\"u}r Angewandte Mathematik und Mechanik},
  14(4):203--212, 1934.

\bibitem{Lindbo2011c}
Dag Lindbo and Anna-Karin Tornberg.
\newblock {Spectral accuracy in fast Ewald-based methods for particle
  simulations}.
\newblock {\em Journal of Computational Physics}, 230(24):8744--8761, oct 2011.

\bibitem{Kolafa1992}
Jiri Kolafa and John~W. Perram.
\newblock {Cutoff Errors in the Ewald Summation Formulae for Point Charge
  Systems}.
\newblock {\em Molecular Simulation}, 9(5):351--368, 1992.

\end{thebibliography}
\bibliographystyle{unsrt}

\section{Appendix}
\label{sec:appendix}

We record the expressions for the first 3 derivatives of $f(r):= \dfrac{\erf(\xi r)}{r}$ that appear in \eqref{eqn:G_R}.
\begin{align*}
& f'(r) = \dfrac{2 \xi \exp(-\xi^2 r^2)}{\sqrt{\pi } r} - \dfrac{\erf(\xi r)}{r^2},\\
& f''(r) = -\dfrac{4 \xi \exp(-\xi^2 r^2)}{\sqrt{\pi } r^2} - \dfrac{4 \xi^3 \exp(-\xi^2 r^2)}{\sqrt{\pi }} + \dfrac{2 \erf(\xi r)}{r^3}, \\
& f'''(r) = \dfrac{4\xi \exp(-\xi^2 x^2)}{\sqrt{\pi}r} \left (\dfrac{3}{r^2} + 2\xi^2 + 2\xi^4r^2  \right) - \dfrac{6\erf(\xi r)}{r^4}.
\end{align*}

Next we record the real and Fourier space parts for the stokeslet, stresslet and rotlet below. Due to the complicated formulae involved, we change our notation slightly and present them directly as quoted by Klinteberg \etal \cite{fsewald2017}. The modification to our case with images and sources is straight-forward and left to the reader.
\begin{subequations}
\begin{align}
  & \stokeslet_{jl}^R(\mb{r}, \xi)   = 2\left( \frac{\xi  e^{-\xi^2 r^2}  }{\sqrt{\pi}} +
      \frac{ \erfc{(\xi r)}  }{2 r} \right) (\delta_{jl} + \hat{r}_j \hat{r}_l) -
    \frac{4\xi}{\sqrt{\pi}}  e^{-\xi^2 r^2}  \delta_{jl} ,  \\
& \mathcal{T}^R_{jlm} (\mb{r}, \xi)  = - \frac{2}{r} \left[ \frac{3 \, \erfc(\xi r)}{r} +
\frac{2 \xi}{\sqrt{\pi}}  \left(  3+2\xi^2r^2\right)
  e^{-\xi^2 r^2} \right]  \hat{r}_j \hat{r}_l  \hat{r}_m +
\frac{4 \xi^3}{\sqrt{\pi}} e^{-\xi^2 r^2} (\delta_{jl}\hat{r}_m + \delta_{lm} \hat{r}_j
+\delta_{mj}\hat{r}_l), \\
 &  W^R_{jl}(\mb{r}, \xi) = 2 \varepsilon_{jlm} \hat{r}_m
\left( \frac{\erfc(\xi r)}{r^2} + \frac{2 \xi}{\sqrt{\pi}}
  \frac{1}{r} e^{-\xi^2 r^2} \right),  
\end{align}
\label{eqn:kernel_real_parts}
\end{subequations}
where $\hat{\mb{r}}=\mb{r}/|\mb{r}|$.

For the Fourier space part we have
\begin{subequations}
\begin{align}
 & \widehat{\stokeslet}^F(\mb{k},\xi)  = \Gwoexp^{\stokeslet}(\mb{k},\xi)
 e^{-k^2/4 \xi^2}\widehat \biharmonic(k), \\ \quad 
 & \widehat{\mathcal{T}}^F(\mb{k},\xi)  = \Gwoexp^{\stresslet}(\mb{k},\xi)
 e^{-k^2/4 \xi^2}\widehat \biharmonic(k), \\ \quad 
 & \widehat{W}^F(\mb{k},\xi)  = \Gwoexp^{\rotlet}(\mb{k},\xi)  e^{-k^2/4 \xi^2}\widehat \harmonic(k), 
\end{align}
\label{eqn:kernel_fourier_space_part}
\end{subequations}
where
\begin{align*}
  \Gwoexp^{\stokeslet}_{jl}(\mb{k},\xi) &= 
  -\left(k^2\delta_{jl} - k_jk_l
  \right) 
  \left(1+k^2/(4\xi^2)\right), \\ 
\Gwoexp^{\stresslet}_{jlm}(\mb{k},\xi) &=
 -i \left[(k_m \delta_{jl} +k_j \delta_{lm} +k_l  \delta_{mj}) k^2   -2
  k_j k_l k_m \right]  
  \left(1+k^2/(4\xi^2)\right) , \\
\Gwoexp^{\rotlet}_{jl} (\mb{k},\xi) &=
2  i \varepsilon_{jlm} k_m. 
\end{align*}

The self-interaction term is non-zero for the stokeslet only, given by
\begin{align*}
\stokeslet^{\mathrm{self}}_{jl} = -\dfrac{4\xi}{\sqrt{\pi}} \delta_{jl}.
\end{align*}


\end{document}